\documentclass[11pt]{article}
\usepackage[T1]{fontenc}
\usepackage[utf8]{inputenc}
\usepackage{lmodern,microtype}
\usepackage[a4paper,margin=27mm]{geometry}
\usepackage{amsmath,amssymb,amsthm,mathtools}
\usepackage{booktabs,longtable,array,enumitem,listings}
\usepackage[hidelinks]{hyperref}
\usepackage{url}
\hypersetup{pdftitle={An exact certificate-based enumeration of Rolle words of degree six},
  pdfauthor={Andrey A. Novikov}}
\newtheorem{theorem}{Theorem}[section]
\newtheorem{lemma}[theorem]{Lemma}
\newtheorem{proposition}[theorem]{Proposition}

\theoremstyle{definition}
\newtheorem{definition}[theorem]{Definition}

\newcommand{\R}{\mathbb R}

\DeclareMathOperator{\spanop}{span}
\usepackage{xcolor}
\setlist{nosep}
\title{An exact certificate-based enumeration\\of Rolle words of degree six}
\author{Andrey A. Novikov}
\date{11 November 2026\\\small Preprint, version 1.1}

\begin{document}
\maketitle
\begin{center}\small\textbf{Preprint --- not peer reviewed.}\\
Most of the mathematical development, code, certificates, and manuscript\\
were generated using ChatGPT (OpenAI), under the direction of the human author.
\end{center}
\begin{abstract}
For a real polynomial of degree $n$ with $n$ distinct real roots, assume
that no two different nonconstant derivatives, including the polynomial
itself, share a root. Recording the derivative orders of all roots from
left to right gives a Rolle word. We present an exact, computer-assisted
classification in degree six: $2706$ words are realizable, with a word
and its reversal counted separately. Starting from $33592$ formal
interlacing words, published necessary restrictions leave $13404$
candidates. Integer realization certificates and $53$ exact exclusion
stages classify every candidate. The surviving set agrees exactly with
the set of certified realizations. All verification uses integer and
rational arithmetic. The accompanying code, certificates, proof notes,
and per-word exclusion ledger make the computation reproducible.
The result concerns degree six; it does not resolve the classification
problem in arbitrary degree or count connected components of coefficient
space. Most of the mathematical development, code, certificates, and
initial manuscript were generated using ChatGPT under the author's
direction. No independent human audit of the complete proof chain has
yet been completed.
\end{abstract}
\noindent\textbf{Keywords:} hyperbolic polynomial; Rolle word; real roots;
derivative; exact enumeration; computer-assisted proof; Bernstein basis.

\section{Problem, context, and result}

Let $P\in\R[x]$ have degree $n$ and $n$ distinct real roots. By Rolle's
theorem, $P^{(j)}$ has $n-j$ distinct real roots for $0\le j<n$.
We impose the additional general-position assumption that roots belonging
to different derivatives are also distinct. Thus the complete collection
contains $N=n(n+1)/2$ distinct real numbers.

\begin{definition}
The \emph{Rolle word} $w(P)$ is the word of length $N$ obtained by sorting
this collection increasingly and replacing a root of $P^{(j)}$ by $j$.
Let $a_n$ denote the number of realizable words. A word and its reversal
are counted separately.
\end{definition}

\paragraph{Relation to earlier work.}
Shapiro and Shapiro~\cite[Section~2]{ShapiroShapiro2012} formulate the
classification problem in this encoding, using the term symbolic
sequence. It also appears as Problem~7 in Shapiro's survey
\cite{Shapiro2015}. A formal word satisfies Rolle interlacing, but need
not be realizable by a polynomial. In degree four, $10$ of the $12$
formal words are realizable. Kostov~\cite[Theorem~2]{Kostov2002} proves
that $a_5=116$, whereas there are $286$ formal words in degree five.
Kostov and Shapiro~\cite{KostovShapiro2002} give further necessary
restrictions. These published results are inputs to the present proof;
in particular, completeness of the quintic classification is not inferred
from a search for examples.

The class of hyperbolic polynomial-like functions is larger than the
class of polynomials considered here. Kostov~\cite{Kostov2007} gives
$236$ realizable arrangements for such functions in degree five,
compared with $116$ for polynomials. The two realization problems must
therefore be distinguished.

There is also a separate asymptotic question, attributed to Arnold by
Shapiro and Shapiro~\cite[Section~2, Problem~1]{ShapiroShapiro2012}:
does the proportion of formal words realizable by polynomials tend to
zero? A recent preprint of Ferudun~\cite{Ferudun2026} addresses that
question by bounding the number of realizable words using
collision resultants and a component bound. It does not give a
classification in degree six. Its asymptotic result is not used below.

\paragraph{Contribution and scope.}
The present work gives a finite classification for degree six, together
with exact realization and exclusion certificates. It provides an
explicit instance of the general classification problem, rather than a
solution in arbitrary degree. The main result is as follows.

\begin{theorem}[Computer-assisted enumeration]\label{thm:main}
Under the definition above, $a_6=2706$.
\end{theorem}

The lower bound comes from $2706$ distinct words with integer
realization certificates. For the upper bound, published restrictions
reduce the formal enumeration to $13404$ candidates, and $53$ checked
exclusion stages remove the other $10698$ words. The final residual set
is compared with the certified witness set, word by word.

The quantity $a_6$ counts total orders of roots, not connected components
of coefficient space. No claim that the realization set of a fixed word
is connected is needed. We make no claim here to have established
absolute publication priority for the count or for each individual
inequality.

\paragraph{Evidence and verification limits.}
The reference computation and the refactored certificate pipeline have
both been executed successfully. They share the specialized mathematical
checkers, so their agreement is a reproducibility check, not an
independent implementation of every proof argument. The lower-bound
certificates are checked directly by integer sign evaluation. The
upper bound additionally depends on the correctness of the necessary
conditions used by the exclusion rules. Successful execution checks the
encoded identities, signs, covers, and set accounting; it does not by
itself validate every mathematical reduction. The development is not
formalized in a proof assistant, and no independent human audit of the
complete proof chain has yet been completed. Section~\ref{sec:ai}
describes the substantial role of ChatGPT and the author's role.

\paragraph{Organization.}
We first enumerate the formal words and apply the published filters,
then describe the realization certificates and the main forms of
exclusion certificate. Sections~\ref{sec:normal} and~\ref{sec:central}
give a common normalization and a complete worked exclusion.
The remaining sections explain the set accounting, reproduction
procedure, and provenance. Appendix~\ref{app:ledger} lists every stage
and its incremental contribution; the ancillary proof notes and source
code supply the stage-specific derivations.

\section{Exhaustive enumeration without root finding}

Write $x^{(j)}_1<\cdots<x^{(j)}_{n-j}$ for the roots of $P^{(j)}$.
The adjacent interlacing relations are
\begin{equation}\label{eq:rolle}
 x^{(j)}_k<x^{(j+1)}_k<x^{(j)}_{k+1}
 \quad(0\le j<n-1,\ 1\le k<n-j).
\end{equation}
Let $c_j$ be the number of occurrences of $j$ in a word prefix.

\begin{lemma}\label{lem:prefix}
A word with multiplicities $n,n-1,\ldots,1$ encodes a linear extension of
\eqref{eq:rolle} if and only if every prefix satisfies
\begin{equation}\label{eq:prefix}
 0\le c_j-c_{j+1}\le1\qquad(0\le j<n-1).
\end{equation}
\end{lemma}
\begin{proof}
After deleting every label except $j$ and $j+1$, strict interlacing says
that the remaining word alternates, beginning and ending with $j$.
This is equivalent to \eqref{eq:prefix} for that pair. Applying this to
every adjacent pair gives exactly all relations \eqref{eq:rolle}.
The occurrence number of a label fixes its root index, so a word cannot
encode two different linear extensions.
\end{proof}

This yields a short depth-first enumerator: append any label whose
multiplicity is not exhausted and whose addition preserves
\eqref{eq:prefix}. There is no numerical feasibility test in this step.
The number of completions can also be computed by the recurrence
\begin{equation}\label{eq:dp}
 T(c)=\sum_{j\text{ admissible at }c}T(c+e_j),
 \qquad T(n,n-1,\ldots,1)=1.
\end{equation}
The implementation memoizes this recurrence and separately generates the
words. Their equality is tested against the archived partially ordered
set enumerator, word by word, for degrees $1$ through $6$.

\begin{center}
\begin{tabular}{lrrrrrr}
\toprule
Degree $n$&1&2&3&4&5&6\\
\midrule
Formal words&1&1&2&12&286&33592\\
Realizable words&1&1&2&10&116&2706\\
\bottomrule
\end{tabular}
\end{center}
The last entry is the result of the full certificate computation, not of
the recurrence alone.

\subsection{Published necessary restrictions}
Removing every zero from a sextic word and subtracting one from the other
labels gives the word of $P'$. The $116$ stored quintic words all have exact
realization certificates. Kostov's theorem \cite[Theorem~2]{Kostov2002}
therefore proves completeness of this particular list. Keeping only
sextic words with one of these tails leaves $13950$ words.

For a degree-$m$ strictly hyperbolic polynomial $f$, Theorem~1 of
\cite{KostovShapiro2002} implies
\begin{equation}\label{eq:ks}
 \left(x^{(0)}_{j+1}<x^{(2)}_j\text{ for }1\le j\le m-2\right)
 \ \Longrightarrow\ x^{(m-3)}_2<x^{(m-1)}_1.
\end{equation}
Apply this to every applicable derivative and to its reflection. Exactly
$13404$ words survive. Both filters are necessary conditions, not
realizability tests.

\section{Exact realizability certificates}

A stored witness consists of an integer polynomial, either given by its
coefficients or by distinct integer roots, a positive integer $q$, and
integers $b_0<\cdots<b_N$. The $i$th letter $j$ asserts that $P^{(j)}$
has a root in $(b_i/q,b_{i+1}/q)$.

\begin{lemma}\label{lem:witness}
If the assigned derivative has opposite nonzero signs at the two
endpoints of every interval, then the witness realizes the asserted word.
\end{lemma}
\begin{proof}
The intermediate value theorem gives a root in every assigned interval.
For $P^{(j)}$ there are exactly $n-j$ such disjoint intervals, equal to
its degree. Therefore it has precisely one root in each, every such root
is simple, and there are no additional roots. Interiors of the intervals
for all derivative orders are mutually disjoint; no root is an endpoint
of its assigned interval. Hence all roots are distinct and have the
specified total order.
\end{proof}

If $f$ has degree $d$, the sign of $f(b/q)$ is the sign of the integer
$q^d f(b/q)$. A denominator-clearing Horner evaluation computes this
integer directly. Thus the verification requires neither approximate
roots nor a tolerance. The witness checker also rejects duplicate
JSON keys, malformed multiplicities, noninteger data, zero denominators,
unordered endpoints, and missing sign changes.

The verified sextic corpus has $2706$ distinct words and is invariant
under reversal. This proves $a_6\ge2706$ independently of the exclusion
calculation.

\section{Soundness of the exclusion certificates}

The exclusion corpus uses several mathematical reductions, but its finite
checks fall into a small number of explicit forms. Necessity of a
reduction is part of the mathematical argument; a valid algebraic identity
cannot repair an invalid necessary inequality.

\subsection{Linear and quadratic identities}
Represent the ordered $21$ nodes by $20$ positive gaps. Root coordinates
relative to the common mean are linear forms in these gaps. Equality of
the root means of successive derivatives gives linear equations $Eg=0$.
Interlacing and additional root inequalities give strict forms $Bg>0$.

\begin{lemma}\label{lem:farkas}
If $y^TE=z^TB$ with $z\ge0$ and $z\ne0$, the strict system
$Eg=0$, $Bg>0$ has no solution.
\end{lemma}
\begin{proof}
The left side evaluated at a putative solution is zero, whereas the right
side is strictly positive.
\end{proof}

The coefficients of these certificates are integers. Quadratic variants
multiply equality constraints by gap coordinates and combine them with
second-moment identities; the positive side is a nonzero nonnegative
combination of products of strict positive forms. Every coefficient of
the resulting polynomial identity is checked exactly.

For completeness, two root inequalities used in these systems are recorded
here. If $x_1<\cdots<x_m$ are the roots of $f$ and
$y_1<\cdots<y_{m-1}$ those of $f'$, then, for $m\ge3$,
\begin{equation}\label{eq:rootcuts}
 \frac{(m-k)x_k+x_{k+1}}{m-k+1}<y_k<
 \frac{x_k+kx_{k+1}}{k+1}\quad(1\le k<m).
\end{equation}
Indeed, split $f'(y_k)/f(y_k)=0$ into terms to the left and right of
$y_k$, then bound each denominator by its nearest endpoint. At least one
bound is strict. In addition, for $1\le k\le m-2$,
\begin{equation}\label{eq:partialsum}
 (k+1)\sum_{i=1}^k y_i<k\sum_{i=1}^{k+1}x_i,
\end{equation}
with the reversed inequality for the corresponding rightmost roots.
To see this, realize the $y_i$ as the eigenvalues of the compression of
$\operatorname{diag}(x_1,\ldots,x_m)$ to the orthogonal complement of the
all-ones vector. The minimum trace on a $k$-dimensional subspace is
$\sum_{i=1}^k y_i$. Restricting to vectors supported on the first $k+1$
coordinates and having zero coordinate sum gives trace
$\frac{k}{k+1}\sum_{i=1}^{k+1}x_i$. Equality would make this proper
subspace invariant under the compression, which is impossible for
distinct $x_i$. The characteristic polynomial identity for this
compression follows directly from the resolvent identity
$f'(z)/f(z)=\sum_i(z-x_i)^{-1}$.

If the roots are centered at their common mean and $S_m$ denotes their
sum of squares, coefficient comparison gives
\begin{equation}\label{eq:moments}
 S_{m-1}=\frac{m-2}{m}S_m,
 \qquad\frac{S_m}{m(m-1)}=\frac{S_2}{2}.
\end{equation}
These are the quadratic equalities used in the initial moment stage.
The same ideas apply to root-subset factors, their derivatives, and
consistent deletion of roots. Metric cuts include rational consequences
of the span bound
\begin{equation}\label{eq:span}
 \spanop(f')\ge\sqrt{\frac{m-2}{m}}\,\spanop(f)
\end{equation}
of Meir--Sharma \cite{MeirSharma1968}. Rational constants are rounded in
the direction that preserves necessity, and the rounding is checked by
integer squares.

\subsection{Polynomial sign identities and Bernstein covers}
Other stages normalize a low-order derivative and integrate it, leaving a
small number of shape parameters and integration constants. A prescribed
word fixes signs of values at marked roots. Eliminating the integration
constants produces necessary polynomial inequalities.

For a polynomial $f$ on a rational box, substitute
$x_i=\ell_i+(r_i-\ell_i)t_i$, $0\le t_i\le1$, and expand in the tensor
Bernstein basis:
\begin{equation}\label{eq:bernstein}
 f(x(t))=\sum_{\alpha\le d} b_\alpha
   \prod_i\binom{d_i}{\alpha_i}t_i^{\alpha_i}(1-t_i)^{d_i-\alpha_i}.
\end{equation}
The basis functions are nonnegative and sum to one. Hence
$f\ge\min_\alpha b_\alpha$ on the entire closed box. All $b_\alpha$ are
rational and are evaluated exactly. Binary subdivision certificates
record an axis and its two halves at each internal node; a leaf identifies
a polynomial with the required sign. The checker verifies coverage as
well as leaf signs. A strict required condition can be contradicted by a
non-strict bound, whereas rejecting an allowed equality boundary requires
a strict bound. The checkers distinguish these cases.

Unbounded regions are handled by explicit changes of variables, positive
coefficient expansions on a tail, or local algebraic identities near a
limiting square. A numerical search failure is never an exclusion
certificate. The archived derivations and the checker source for every
stage are indexed by the public registry.

\subsection{Hermite constraints}
For a monic real polynomial $f$ of degree $m$, let $s_k$ denote its Newton
sums, with $s_0=m$. The Hermite matrix is $H_m=(s_{i+j})_{0\le i,j<m}$.
For distinct real roots it is a positive definite Vandermonde Gram matrix.
Conversely, if a pair of roots is nonreal, real-coefficient interpolation
gives a polynomial taking values $\mathrm i,-\mathrm i$ on that pair and
zero at the other distinct roots. Its Hermite quadratic form is negative.
Repeated roots make the matrix singular. This proves the strict
hyperbolicity criterion used by the Hermite stages. Newton identities,
matrix test-vector identities, and sign bounds for their coefficients are
checked in exact arithmetic. Low-order matrix tests are only necessary
conditions; passing them does not certify a sextic witness.

\section{A common normalization for the final exclusions}\label{sec:normal}

Multiplication of $P$ by $-1$ preserves its word, so we may assume a
positive leading coefficient. Normalize the two middle roots of $P''$ to
$0$ and $1$, and write the outer roots as $-u$ and $T$, where $u>0$ and
$T>1$. After a positive scaling of $P$, write
\begin{align}
 R=P''&=(x+u)x(x-1)(x-T), &T&=1+u+2r,\\
 H&=u^2+u+2ur-2r-1, &P&=F-qx+C,\\
 F(x)&=\frac{x^6}{30}-\frac{(1+r)x^5}{10}
       -\frac{Hx^4}{12}+\frac{uTx^3}{6}, &Q&=F'.
\end{align}
Let $a<c<d$ be the roots of $R'$, $v_-<v_+$ those of $R''$, and
$\mu=(1+r)/2$ the root of $R'''$. Reflection permits the orientation
$c<\mu$, equivalent to $r>0$, since
\begin{equation}
 R'(\mu)=-r(r+u)(r+u+1).
\end{equation}
General position excludes equality. A strictly hyperbolic sextic must
have negative first critical minimum and positive fourth critical
maximum. These elementary sign requirements are used repeatedly.

The first $46$ stages leave $2734$ candidates. The last seven stages
exclude $28$ of these candidates, leaving the certified witness set.
The first two share a recorded recovery stage; their
net contributions are nevertheless counted separately.
\begin{center}\small
\begin{tabular}{p{.55\linewidth}rr}
\toprule
Rule&New words&Cover leaves\\
\midrule
Negative left cubic value&2&7\\
Matched left-minimum/right-maximum regimes&2&$241+530$\\
Negative right maximum&2&34\\
Positive return at the mean&4&18\\
Negative inflection values&6&17\\
First minimum with $H>0$&6&81\\
Central dip and first minimum&6&194\\
\bottomrule
\end{tabular}
\end{center}
Leaf counts describe finite covers, not the number of independent
mathematical hypotheses. Tail inequalities and limiting identities are
additional checked parts of the relevant rules. We give the final rule
in detail, to exhibit the link between a word and its finite certificate.

\section{The final central-dip obstruction}\label{sec:central}

\begin{proposition}\label{prop:central}
No strictly hyperbolic sextic in the normalization of
Section~\ref{sec:normal} can satisfy all the conditions
\begin{equation}\label{eq:central-signs}
 \begin{gathered}
 R(v_-)>0,\quad R(\mu)>0,\quad R(v_+)<0,\\
 P'(a)<0,\quad P(c)<0,\quad P'(v_+)>0,\quad P(0)>0.
 \end{gathered}
\end{equation}
\end{proposition}

\begin{proof}
Interlacing and $R''(0)=-2H$ imply
\begin{equation}\label{eq:central-domain}
 H<0,\qquad 0<r<1,\qquad v_+>1,\qquad H+6r>0.
\end{equation}
The last implication uses $R''(1)=-12r-2H<0$.
Since $H$ is increasing in $u$ and $H(1,r)=1$, we have $u<1$.
If $s_1<-u$ is the first root of $P'$, the necessary inequalities are
\begin{equation}\label{eq:q-constraints}
 \begin{gathered}
 C>0,\quad P(s_1)<0,\quad 0<q<\min\{Q(-u),Q(1)\},\quad q=Q(s_1),\\
 q>Q(a),\qquad q>F(c)/c,\qquad q<Q(v_+).
 \end{gathered}
\end{equation}
The middle strict lower bound follows from $P(c)<0$ and $C>0$.

Direct evaluation gives
\begin{align}
 R'(-3u/4)&=-\frac{u(6ru+16r+3u^2+14u+8)}{16}<0,\\
 R'(-u/2)&=\frac{u^2(r+u+1)}{2}>0,\\
 Q(-3u/2)&=-\frac{9u^4(5r+7u+5)}{160}<0.
\end{align}
Consequently $a=-ut$, $1/2<t<3/4$, and $s_1=-uv$, $1<v<3/2$.
Furthermore $R'(1/2)=r/2>0$, so $c>1/2$. On $(0,1)$ the function $Q$
increases strictly, and $F(0)=Q(0)=0$; hence $F(x)/x$ increases strictly.
Thus $q>2F(1/2)$, whereas $q<Q(-u)$. On
$0\le u\le2/5$, $0\le r\le1$, the minimum Bernstein coefficient of
$2F(1/2)-Q(-u)$ is $2143/1000000>0$. This exact calculation excludes
that rectangle and gives $u>2/5$.

Define the following polynomials in $u,t$:
\begin{align}
 D&=ut(3t-2)+2t-1,\\
 N&=(2ut+1)\{u(1-2t^2)-2t+1\},\\
 D_N&=2D+N, &T_N&=(1+u)D+N,\\
 H_N&=D(u^2+u-1)+N(u-1), &H_b&=H_N/u.
\end{align}
The division defining $H_b$ is exact in the polynomial ring. The
logarithmic derivative of $R$ at $a=-ut$ gives
\begin{equation}
 D=\frac{ut(1-t)(ut+1)}{ut+T}>0,
 \qquad r=\frac{N}{2D},\qquad H_N=DH.
\end{equation}
It follows that
\begin{equation}\label{eq:necessary-D}
 D>0,\quad N>0,\quad N<2D,\quad H_b<0,\quad H_N+3N>0.
\end{equation}
The full necessary parameter region is contained in the closed box
\begin{equation}\label{eq:central-box}
 [2/5,1]\times[1/2,3/4]\times[1,3/2]
 \quad\text{in }(u,t,v).
\end{equation}

Set
\begin{align}
 E(x)&=60DF(x)=2Dx^6-3D_Nx^5-5H_Nx^4+10uT_Nx^3,\\
 Q_A&=E'(-ut)/u^4,\qquad Q_X=E'(-uv)/u^4,\\
 G_q&=\{E(-uv)+uvE'(-uv)\}/u^5.
\end{align}
All displayed quotients are polynomials. We must have
$Q_X>Q_A$, $Q_X>0$, and $G_q<0$, the last because
$P(s_1)=u^5G_q/(60D)+C<0$. Nonnegativity of any polynomial in
\begin{equation}\label{eq:central-base-polys}
 -D,\ -N,\ N-2D,\ H_b,\ -(H_N+3N),\ Q_A-Q_X,\ -Q_X,\ G_q
\end{equation}
therefore yields a contradiction.

Two further alternatives use certified rational root bounds. A constant
$c_L\in[1/2,1)$ satisfying $2DR'(c_L)\ge0$ uniformly on a parameter box
obeys $c_L\le c$; $c_L=1/2$ is available without a separate evaluation.
Hence nonnegativity of
\begin{equation}\label{eq:central-cL}
 E(c_L)/c_L-E'(-uv)
\end{equation}
contradicts $q>F(c)/c\ge F(c_L)/c_L$.
Similarly a constant $w\in[1,2)$ satisfying $-2DR''(w)\ge0$ obeys
$w\le v_+$; the initial bound $w=1$ follows from
\eqref{eq:central-domain}. Since $Q$ decreases on $[1,v_+]$, nonnegativity of
\begin{equation}\label{eq:central-w}
 E'(-uv)-E'(w)
\end{equation}
contradicts $q<Q(v_+)\le Q(w)$.

The attached certificate partitions \eqref{eq:central-box} into $194$
rational boxes, using $387$ tree nodes and maximum depth $28$. For the
ten alternatives in \eqref{eq:central-base-polys},
\eqref{eq:central-cL}, and \eqref{eq:central-w}, respectively, the leaf
counts are
\begin{equation}
 (0,3,3,5,3,63,0,72,22,23).
\end{equation}
The checker reconstructs the polynomials from $E$, verifies the
rationalization identities, checks each auxiliary root bound, and
checks every Bernstein coefficient with rational arithmetic. All of
\eqref{eq:central-box}, including its boundary, is covered. Every possible
parameter point therefore violates a necessary strict inequality.
\end{proof}

The six final words satisfying \eqref{eq:central-signs} are
\begin{center}\ttfamily\small
\begin{tabular}{ll}
010213204130524103210&012301425031402312010\\
010213204130524130210&012031425031402312010\\
010213204130524132010&010231425031402312010
\end{tabular}
\end{center}
Each row is a reflection pair. The predicate is independently checked
by reading derivative signs directly from word positions. It matches
$390$ initial candidates, of which $384$ had already been excluded;
its net contribution is exactly six.

\section{Set accounting and the proof of the count}

Let $C_0$ be the $13404$ words left after the published filters. For
stage $i$, let $M_i$ be the set recognized by its checked obstruction.
The actual exclusion and residual sets are defined by
\begin{equation}\label{eq:ledger}
 E_i=C_{i-1}\cap M_i,\qquad C_i=C_{i-1}\setminus E_i.
\end{equation}
The sets $E_i$ are disjoint by construction. A rule can match many words
already excluded without changing its contribution. The program checks
that no $M_i$ intersects the verified witness set, and that the final
residual set equals that witness set exactly, not merely in cardinality.
For historical predicate stages the individual exclusion lists are also
compared with their archived counterparts. Expected counts serve only as
regression checks and are never used to select excluded words.

\begin{proof}[Proof of Theorem~\ref{thm:main}]
Lemma~\ref{lem:prefix} enumerates every formal word exactly once.
The quintic theorem and \eqref{eq:ks} leave $C_0$, which contains every
realizable sextic word. Each checked exclusion rule is a necessary-sign
contradiction of the kinds described above, so no realizable word is
removed by \eqref{eq:ledger}. The $53$ stages remove $10698$ distinct
words and leave precisely the $2706$ words with integer realization
certificates. Lemma~\ref{lem:witness} proves all these realizations.
Thus both bounds are $2706$.
\end{proof}

The full incremental ledger is given in Appendix~\ref{app:ledger}.
After stage $46$, the residual set has $2734$ words. The last seven
stages contribute $2,2,2,4,6,6,6$ new exclusions. Some historical proof
notes retain the intermediate bound $2734$; the complete ledger gives
the final count $2706$.

\section{Software structure and reproduction}

The readable front end is the package \texttt{rolle6}. Its responsibilities
are separated as follows:
\begin{center}\small
\begin{tabular}{lp{.64\linewidth}}
\toprule
Module&Responsibility\\
\midrule
\texttt{words.py}&Prefix enumeration, dynamic-programming counts, tail extraction, and published filters.\\
\texttt{witnesses.py}&Integer evaluation and exact interval witness validation.\\
\texttt{stages.py}&Explicit registry of $53$ proof stages and their reference increments.\\
\texttt{backend.py}&Invocation of the archived exact proof checkers, with a per-run cache.\\
\texttt{pipeline.py}&Proof-before-exclusion control flow, set ledger, input hashes, and machine-readable results.\\
\bottomrule
\end{tabular}
\end{center}
The specialized mathematical modules in \texttt{research/} remain unchanged.
Their hashes are compared with the completed pre-refactoring archive.
This preserves the proof corpus while replacing the enumeration and
monolithic bookkeeping. It is not a claim that every algebraic checker
has been independently reimplemented.

From the root of the ancillary code directory, use Python~3.10 or later:
\begin{lstlisting}
python3 -m rolle6 enumerate --degree 6
python3 -m rolle6 verify --output-dir verification_output
python3 -m unittest discover -s tests -v
\end{lstlisting}
No third-party Python package is required for verification. A complete
run writes the proof log, each stage's exact exclusion list, the final
realizable list, and a dependency hash inventory. The verifier refuses
optimized Python execution, because the archived checkers use assertions
that \texttt{python -O} would disable.

The refactoring tests compare every generated word with the earlier
partially ordered set enumerator in degrees $1$ through $6$, compare every
published-filter decision, check all $2822$ quintic and sextic witnesses,
and check rejection of deliberately damaged certificates. The entire
exclusion pipeline is then run against the archived data. Search scripts,
random seeds, floating-point optimizers, and intermediate unsuccessful
experiments are not required to reproduce the result.

\section{Substantial use of ChatGPT and verification limits}\label{sec:ai}

\textbf{Most of the mathematical development and computational work in this
project was generated by ChatGPT (OpenAI), not merely edited or translated
with it.} This includes proposing many of the reductions, drafting most
proof arguments, writing and revising the search and verification
programs, constructing exact certificates, analyzing the remaining cases,
performing the code refactoring, and drafting most of this manuscript.
This disclosure concerns the mathematical and computational development
as well as the writing.

The human author selected the research problem, directed the interactive
investigation, proposed and questioned lines of attack, requested checks
against overlap and missed cases, required intermediate results and checks
to be recorded, and requested explicit disclosure of the AI contribution.
The project was carried out through ChatGPT during November 2026. Exact
model-version attribution is not available for every part of the extended
session; no unsupported per-step model identifier or contribution
percentage is asserted.

ChatGPT is disclosed as a research and writing tool. The human author
remains responsible for any submitted version. In particular, repeated
checks produced within the same AI-assisted workflow are not a substitute
for an independent mathematical and software audit. At the time of this
preprint, no such independent human audit of the entire chain has been
completed. Publishing the exact certificates, explicit assumptions,
readable enumeration, per-word ledger, and unchanged verification backend
is intended to make that audit feasible.

Execution of the exact checker, correctness of the mathematical
reductions, and novelty are separate questions. The successful runs
establish reproducibility of the encoded computation. The mathematical
arguments and published inputs support the reductions. The literature
comparison in Section~1 places the result in context but does not
constitute a certification of publication priority.
The manuscript does not attribute classical results, the known quintic
classification, or published root inequalities to the AI-assisted project.

\appendix
\section{Complete incremental exclusion ledger}\label{app:ledger}
The initial candidate count is $13404$. Each row lists only newly removed
words. The machine-readable registry supplies the corresponding checker
and explanatory source; the verification report supplies the words themselves.
\small
\begin{longtable}{r p{.61\textwidth} r r}
\toprule
No. & Rule & New & Remaining\\
\midrule
\endfirsthead
\toprule
No. & Rule (continued) & New & Remaining\\
\midrule
\endhead
\bottomrule
\endfoot
1 & Linear root inequalities & 2718 & 10686 \\
2 & Quadratic root moments & 2282 & 8404 \\
3 & Products of positive linear forms & 4 & 8400 \\
4 & Quintic primitive types & 1168 & 7232 \\
5 & Primitive sign at the mean & 396 & 6836 \\
6 & First and last interval averages & 534 & 6302 \\
7 & Middle interval and convexity & 558 & 5744 \\
8 & Central-value identities & 298 & 5446 \\
9 & Parametric middle-root identities & 212 & 5234 \\
10 & Polar derivative at the mean & 12 & 5222 \\
11 & Deletion of one root & 44 & 5178 \\
12 & Joint deletion of roots & 8 & 5170 \\
13 & Subfactor sign constraints & 16 & 5154 \\
14 & Metric span constraints & 6 & 5148 \\
15 & Coupled second moments & 2 & 5146 \\
16 & Quartic middle-root bridge & 82 & 5064 \\
17 & Outer cubic roots & 6 & 5058 \\
18 & Two-parameter polynomial bounds & 6 & 5052 \\
19 & Right critical maximum & 12 & 5040 \\
20 & Stationary-point patterns & 1152 & 3888 \\
21 & Inflection-point patterns & 790 & 3098 \\
22 & Full stationary constraints & 30 & 3068 \\
23 & Covering boundary charts & 16 & 3052 \\
24 & Left minimum and inflection & 2 & 3050 \\
25 & Coupled root constraints & 142 & 2908 \\
26 & Pairs of inflection points & 38 & 2870 \\
27 & Hermite fibres & 8 & 2862 \\
28 & Hermite endpoint bounds & 10 & 2852 \\
29 & Laguerre integral identities & 6 & 2846 \\
30 & Outer secants & 44 & 2802 \\
31 & Marker quadrature & 2 & 2800 \\
32 & Hermite marker matrices & 2 & 2798 \\
33 & Order-six Hermite markers & 4 & 2794 \\
34 & Marker boundary charts & 14 & 2780 \\
35 & Moving Hermite markers & 4 & 2776 \\
36 & Quintic minimum envelope & 16 & 2760 \\
37 & Three inflection constraints & 2 & 2758 \\
38 & Middle triple identity & 4 & 2754 \\
39 & Right fourth-derivative minimum & 2 & 2752 \\
40 & Outer inflection average & 2 & 2750 \\
41 & Central cubic chord & 4 & 2746 \\
42 & Left cubic chord & 2 & 2744 \\
43 & Left cubic and mean average & 4 & 2740 \\
44 & Fourth-derivative average & 2 & 2738 \\
45 & Fourth-derivative chord & 2 & 2736 \\
46 & Outer cubic chord & 2 & 2734 \\
47 & Negative left cubic value & 2 & 2732 \\
48 & Matched extrema regimes & 2 & 2730 \\
49 & Negative right maximum & 2 & 2728 \\
50 & Positive return at the mean & 4 & 2724 \\
51 & Three negative inflection values & 6 & 2718 \\
52 & First minimum with positive H & 6 & 2712 \\
53 & Central dip and first minimum & 6 & 2706 \\
\end{longtable}
\normalsize

\section{Minimal enumeration pseudocode}
The following equivalent form emphasizes the prefix invariant. The actual
implementation also checks input types and supports a memoized counting mode.
\begin{lstlisting}[language=Python]
def enumerate_words(n):
    target = tuple(range(n, 0, -1))

    def visit(counts, prefix):
        if counts == target:
            yield prefix
            return
        for j in range(n):
            if counts[j] == target[j]:
                continue
            updated = list(counts)
            updated[j] += 1
            if all(0 <= updated[k] - updated[k+1] <= 1
                   for k in range(n-1)):
                yield from visit(tuple(updated), prefix + str(j))

    yield from visit((0,) * n, "")
\end{lstlisting}

\clearpage
\bibliographystyle{plain}
\bibliography{main}

@article{Kostov2002,
  author = {Kostov, Vladimir Petrov},
  title = {Discriminant sets of families of hyperbolic polynomials of degree 4 and 5},
  journal = {Serdica Mathematical Journal},
  year = {2002}, volume = {28}, number = {2}, pages = {117--152},
  url = {https://www.math.bas.bg/serdica/2002/2002-117-152.pdf},
  note = {\url{https://www.math.bas.bg/serdica/2002/2002-117-152.pdf}}
}

@article{KostovShapiro2002,
  author = {Kostov, Vladimir Petrov and Shapiro, Boris Zalmanovich},
  title = {On arrangements of roots for a real hyperbolic polynomial and its derivatives},
  journal = {Bulletin des Sciences Math\'ematiques},
  year = {2002}, volume = {126}, number = {1}, pages = {45--60},
  doi = {10.1016/S0007-4497(01)01106-X},
  url = {https://staff.math.su.se/shapiro/Articles/kostov.pdf},
  note = {\href{https://doi.org/10.1016/S0007-4497(01)01106-X}{doi:10.1016/S0007-4497(01)01106-X}}
}

@article{MeirSharma1968,
  author = {Meir, A. and Sharma, A.},
  title = {On zeros of derivatives of polynomials},
  journal = {Canadian Mathematical Bulletin},
  year = {1968}, volume = {11}, number = {3}, pages = {443--445},
  doi = {10.4153/CMB-1968-051-x},
  note = {\href{https://doi.org/10.4153/CMB-1968-051-x}{doi:10.4153/CMB-1968-051-x}}
}

@article{ShapiroShapiro2012,
  author = {Shapiro, Boris and Shapiro, Michael},
  title = {A Few Riddles Behind {Rolle}'s Theorem},
  journal = {The American Mathematical Monthly},
  year = {2012}, volume = {119}, number = {9}, pages = {787--793},
  doi = {10.4169/amer.math.monthly.119.09.787},
  note = {\href{https://doi.org/10.4169/amer.math.monthly.119.09.787}{doi:10.4169/amer.math.monthly.119.09.787}}
}

@article{Shapiro2015,
  author = {Shapiro, Boris},
  title = {Problems Around Polynomials: The Good, The Bad and The Ugly\ldots},
  journal = {Arnold Mathematical Journal},
  year = {2015}, volume = {1}, number = {1}, pages = {91--99},
  doi = {10.1007/s40598-015-0008-4},
  note = {\href{https://doi.org/10.1007/s40598-015-0008-4}{doi:10.1007/s40598-015-0008-4}}
}

@article{Kostov2007,
  author = {Kostov, Vladimir Petrov},
  title = {On root arrangements for hyperbolic polynomial-like functions and their derivatives},
  journal = {Bulletin des Sciences Math\'ematiques},
  year = {2007}, volume = {131}, number = {5}, pages = {477--492},
  doi = {10.1016/j.bulsci.2006.12.004},
  note = {\href{https://doi.org/10.1016/j.bulsci.2006.12.004}{doi:10.1016/j.bulsci.2006.12.004}}
}

@misc{Ferudun2026,
  author = {Ferudun, Alper},
  title = {Vanishing Density of Realizable {Rolle} Words},
  year = {2026},
  howpublished = {Zenodo preprint, version 1.0, 25 October 2026},
  doi = {10.5281/zenodo.23647439},
  note = {\href{https://doi.org/10.5281/zenodo.23647439}{doi:10.5281/zenodo.23647439}}
}
\end{document}